\documentclass[oneside,10pt]{amsart}
\usepackage{amsmath}
\usepackage{amssymb}
\usepackage{amsthm}
\usepackage{float}
\usepackage{ulem}
\usepackage{xcolor}
\usepackage[utf8x]{inputenc}
\usepackage{comment}
\usepackage{framed}
\usepackage{caption}
\DeclareMathAlphabet{\mathbbb}{U}{bbold}{m}{n}

\newtheorem{thm}{Theorem}[section]
\newtheorem{cor}[thm]{Corollary}

\theoremstyle{definition}
\newtheorem{defn}[thm]{Definition}
\theoremstyle{remark}
\newtheorem{rem}[thm]{Remark}

\numberwithin{equation}{section}
\usepackage{graphicx}

\usepackage{enumerate}

\newcommand{\matrice}{\begin{pmatrix}}
\newcommand{\ok}{\end{pmatrix}}

\begin{document}
\title[]{Semiclassical eigenvalue bounds for compact homogeneous irreducible Riemannian manifolds}
\thanks{}

\author[Provenzano]{Luigi Provenzano}
\address{Dipartimento di Scienze di Base e Applicate per l'Ingegneria, Sapienza Universit\`a di Roma, Via Scarpa 12 - 00161 Roma, Italy, e-mail: {\sf luigi.provenzano@uniroma1.it}.}
\author[Stubbe]{Joachim Stubbe}
\address{School of Basic Sciences, EPFL, Rte Cantonale, 1015 Lausanne, Switzerland, e-mail: {\sf joachim.stubbe@epfl.ch}.}

\begin{abstract}
We exploit an identity for the gradients of Laplacian eigenfunctions on compact homogeneous Riemannian manifolds with irreducible linear isotropy group to obtain asymptotically sharp universal eigenvalue inequalities and sharp Weyl bounds on Riesz means. The approach is non variational and is based on identities for spectral quantities in the form of sum rules.
\end{abstract}

\keywords{Homogeneous space, irreducible isotropy group, sharp eigenvalue bounds, universal inequalities, sum rules}
\subjclass{58C40, 58J50, 35P15, 35P20}

\maketitle

%\section{An identity on tangent spaces to homogeneous Riemannian manifolds}

\section{Introduction and statement of the main results}
Let $\mathcal{M}$ be a $d$-dimensional, compact homogeneous Riemannian manifold, and let $\Delta$ be the (negative) Laplacian on  $\mathcal{M}$. The spectrum of  $-\Delta$ consists of non-negative {\it eigenvalue levels}
\begin{equation}\label{eigenvalue-levels}
 0=\Lambda_0<\Lambda_1<\cdots<\Lambda_l<\cdots\nearrow \infty
\end{equation}
with multiplicities
\begin{equation}\label{multiplicities}
 1=M_0<M_1\leq\cdots\leq M_l\leq\cdots
\end{equation}
An orthonormal basis of the corresponding eigenspaces consists of $M_l$  eigenfunctions in $L^2(\mathcal{M})$ which we denote by $Y_l^k$, $k=1,\ldots  M_l$.  We may assume that  $Y_l^k$ are real-valued. 
In \cite{G75} the following addition formula was shown: for $x,y\in \mathcal{M}$ the function
\begin{equation}\label{addition-formula1}
 Z_l(x,y)=\sum_{k=1}^{M_l}Y_l^k(x)Y_l^k(y)
\end{equation}
is a zonal eigenfunction of the Laplacian in the corresponding eigenspace. In particular,
\begin{equation}\label{addition-formula2}
 Z_l(x,x)=\sum_{k=1}^{M_l}|Y_l^k(x)|^2=\frac{M_l}{|\mathcal{M}|},
\end{equation}
where $|\mathcal M|$ is the volume of $\mathcal M$. The second identity in \eqref{addition-formula2} expresses the homogeneity of the manifold on the spectral level: the density of a Laplacian multiplet is constant over the manifold. Formula \eqref{addition-formula2} has several consequences. For example, taking the Laplacian on both sides of \eqref{addition-formula2} we get
\begin{equation}\label{addition-formula3}
 \sum_{k=1}^{M_l}|\nabla Y_l^k(x)|^2=\frac{M_l\Lambda_l}{|\mathcal{M}|}
\end{equation}
where $\nabla$ denotes the gradient operator and $|V|$ the length of the vector $V$ in the tangent space at a point $x\in\mathcal{M}$.

In the present note we consider compact homogeneous Riemannian manifolds for which the linear isotropy group is irreducible on the tangent space (see Definition \ref{def_iso_irr}) and investigate the effect of this property at a spectral level. If $\mathcal M$ is a compact, rank one symmetric space (CROSS), it falls in this class. We recall that CROSS are all the compact {\it isotropic} spaces, namely the spaces for which the stabilizer of $x$ (the subgroup of isometries fixing $x$)  acts transitively on $T_x\mathcal M$, see e.g., \cite[\S 8.12]{wolf_book}. While homogeneity means that the manifold looks the same at each point, isotropy means that it looks the same in any direction (at any point). Hence isotropy implies homogeneity, which is a well-known result.

The manifolds that we consider, however, satisfy a weaker assumption than isotropy (which is quite  restrictive). Roughly speaking, for a $d$-dimensional compact homogeneous Riemannian manifold with irreducible isotropy group, given any direction, there exist at least other $d-1$ linearly independent directions along which the space looks the same. This vague statement can be clarified by means of the simple example of the $2$-dimensional flat torus for which we present some explicit computations in  Appendix \ref{AppB}. In fact, any flat torus is homogeneous but not isotropic. Nevertheless, the square torus and the equilateral torus are isotropy irreducible (they are the only such tori).

The crucial observation leading to our main result (Theorem \ref{main_thm_2}) is the following ``addition-type formula'', or identity on tangent spaces, which holds for compact homogeneous spaces with irreducible linear isotropy group.
\begin{thm}\label{main_thm}
Let $\mathcal{M}$ be an isotropy irreducible compact Riemannian manifold of dimension $d$.  Let $x\in\mathcal{M} $ and $V \in T_x\mathcal{M}$. If we denote by $\langle\cdot,\cdot\rangle$ the scalar product in the tangent space $T_x\mathcal{M}$ induced by the Riemannian metric on $\mathcal M$ then
\begin{equation}\label{addition-formula4}
\sum_{i=1}^k\langle V,\nabla Y_l^k(x)\rangle\nabla Y_l^k(x)=\frac{M_l\Lambda_l}{d|\mathcal M|}V.
\end{equation}
In particular,
  \begin{equation}\label{addition-formula5}
\sum_{k=1}^{M_l}|\langle V,\nabla Y_l^k(x)\rangle|^2=\frac{M_l\Lambda_l}{d|\mathcal{M}|}\, |V|^2.
\end{equation}
\end{thm}

Let us briefly comment on identity \eqref{addition-formula4}. It states that the gradients of the eigenfunctions of any orthonormal basis of an eigenspace form a {\it tight frame} for $T_x\mathcal M$. We remark that identity \eqref{addition-formula4} is related with the theory of minimal isometric immersions of compact manifolds into spheres using eigenfunctions as coordinates, see e.g., \cite{takahashi} (see also \cite{Li81}). 

It is worth noting that an analogous identity holds in $\mathbb R^d$, even though the spectrum is not discrete. In fact, the spectrum of the Laplacian in $\mathbb R^d$ is $[0,+\infty)$, and for $\lambda>0$ a family of generalized eigenfunctions is given by $u_{\lambda,\xi}(x)=e^{i\sqrt{\lambda}\langle \xi,x\rangle}$, with $\xi\in\mathbb R^d$, $|\xi|=1$. We have then
$$
\int_{\mathbb S^{d-1}}|f_{\lambda,\xi}(x)|^2\,d\xi=|\mathbb S^{d-1}|,
$$
for all $x\in\mathbb R^n$ (this is reminiscent of the classic addition formulas, in particular, of \eqref{addition-formula2}). Now, for a vector $V\in\mathbb R^d$ we compute
$$
\int_{\mathbb S^{d-1}}\langle V,\nabla f_{\lambda,\xi}\rangle \overline{\nabla  f}_{\lambda,\xi}\,d\xi=\frac{|\mathbb S^{d-1}|\lambda}{d}V,
$$
which is the analogous of \eqref{addition-formula4}. We include a proof of Theorem \ref{main_thm} in Appendix \ref{AppA}.

The use of Theorem \ref{main_thm} was motivated by the trace identities for eigenvalues of the Laplacian or more general operators $H$ with purely discrete spectrum derived in \cite{HaSt0,HaSt11} which the authors called sum rules refering to early quantum mechanics (see e.g., \cite{BJ,HaSt0}). The basic idea of \cite{HaSt0,HaSt11} is that of exploiting algebraic relations among the first and second commutators of
$H$ with an auxiliary self-adjoint operator $G$, which is usually a multiplication operator by a suitable function (which by abuse of notation we still denote by $G$). Concretely, if $H$ is the Dirichlet Laplacian on a domain $\Omega$ of $\mathbb R^d$, exploiting the trace identities in \cite{HaSt11} one ends up with terms of the form $\displaystyle \int_{\Omega}|\langle\nabla \phi,\nabla G\rangle|^2\,dx$,  where $\phi$ is a Dirichlet Laplacian eigenfunction on $\Omega$. One then chooses $G=x_j$, where $x_j$, $j=1,...,d$ are the  Cartesian coordinate functions, and averages over $j$. As we shall see below, in our situation the natural choice (analogous to the Euclidean case) will be to take $G=Y_l^k$ as coordinates and to average over the multiplet $Y_l^k$, $k=1,..., M_l$. At this point formula \eqref{addition-formula5} will be crucial.  The trace identities found in \cite{HaSt0} lead to quadratic inequalities of the form of \eqref{Polynomial-inequality-Dirichlet-Laplacian} here below.  In turn, these kind of quadratic inequalities imply asymptotically sharp universal eigenvalue bounds and, by the discovery made in \cite{HaHe0},  sharp Weyl-type bounds for Riesz-means or  sharp Lieb-Thirring inequalities \cite{St0} for Schr\"{o}dinger operators. Sum rules for manifolds immersed in the Euclidean space were first derived in \cite{EHI,Ha0}  and later in \cite{HaSt11}, where also an algebraic identity for the spectrum of an abstract selfadjoint operator $H$ defined on a Hilbert space $\mathcal{H}$ was shown. In \cite{EHI} the authors discuss universal inequalities derived from sum rules on compact manifolds isometrically immersed in $\mathbb R^n$. In the present note we derive universal inequalities for compact homogeneous irreducible Riemannian manifolds, which can be isometrically immersed in $\mathbb S^n$ by eigenfunctions, and for their domains.

We can now discuss the main application of Theorem \ref{main_thm}. Let  $\Omega$ be a domain in $\mathcal M$ and let $\lambda_k=\lambda_k(\Omega)$ be the $k$-th eigenvalue of the Dirichlet Laplacian on  $\Omega$. When $\Omega=\mathcal{M} $ (and we have no boundary) the eigenvalues $\lambda_k(\mathcal M)$ are just the Laplace eigenvalues on $\mathcal{M} $ and we have, according to our previous notation,

\begin{equation}\label{eigenvalues-on-M}
 0=\lambda_1(\mathcal M)<\lambda_2(\mathcal M)=\cdots=\lambda_{M_1+1}(\mathcal M)<\lambda_{M_1+2}(\mathcal M)=\cdots
\end{equation}
with $\lambda_2(\mathcal M)=\Lambda_1$, $\lambda_{M_1+1}(\mathcal M)=\Lambda_{2}$, etc. In the following we shall omit the dependence of $\lambda_k$ on $\Omega$ when unnecessary.

 We define the quadratic polynomials $P_N(z)$ and $Q_N(z)$ by
\begin{equation}\label{P-N-of-z-Dirichlet-Laplacian}
\begin{split}
    P_N(z) & =\sum_{j=1}^{N}(z-\lambda_j)(z-\Lambda_1-\frac{d+4}{d}\,\lambda_j) \\
     & =Nz^2-2\,\frac{d+2}{d}\big(\sum_{j=1}^{N}\lambda_j\big)z-\Lambda_1Nz+ \frac{d+4}{d}\,\sum_{j=1}^{N}\lambda^2_j +\Lambda_1\sum_{j=1}^{N}\lambda_j,\\
\end{split}
\end{equation}
\begin{equation}\label{Q-N-of-z-Dirichlet-Laplacian}
  Q_N(z)=N(z-\lambda_N)(z-\lambda_{N+1}).
\end{equation}

We have the following Theorem:

\begin{thm}\label{main_thm_2}
Let $\mathcal{M}$ be an isotropy irreducible compact Riemannian manifold of dimension $d$.  Let $\Omega\subset \mathcal{M} $ and $\lambda_k=\lambda_k(\Omega)$ be the $k$-th eigenvalue of the Dirichlet Laplacian on  $\Omega$.
 Then, for all $z\in [\lambda_N,\lambda_{N+1}]$ the inequality
\begin{equation}\label{Polynomial-inequality-Dirichlet-Laplacian}
  P_N(z)\leq Q_N(z).
\end{equation} 
holds.
\end{thm}

Theorem \ref{main_thm_2} is proved in Section \ref{proof1}. As already said, the quadratic inequality \eqref{Polynomial-inequality-Dirichlet-Laplacian} has several consequences as shown e.g., in \cite{HaSt11}: universal inequalities for eigenvalue gaps and for sums of eigenvalues and its squares, as well as for Riesz-means $R_{\sigma}(z)$ which are defined by
\begin{equation}\label{Riesz-means-Dirichlet-Laplacian}
 R_{\sigma}(z)= \sum_{k} (z-\lambda_k)_{+}^{\sigma}
\end{equation}
for $\sigma\geq 0$ where $a_{+}$ denotes the positive part of the real number $a$. In particular, we have the following Corollary (see \cite{HaSt11}):

\begin{cor}\label{cor-main_thm_2}
Let $\mathcal{M}$ be an isotropy irreducible compact Riemannian manifold of dimension $d$.  Let $\Omega\subset \mathcal{M} $ and $\lambda_k=\lambda_k(\Omega)$ be the $k$-th eigenvalue of the Dirichlet Laplacian on  $\Omega$. 
The corresponding Riesz-mean $ R_{2}(z)$ satisfies the differential inequality
\begin{equation}\label{R-2-diff-inequality}
  \frac{d}{dz}\frac{R_2(z)}{(z+\frac{d}{4}\,\Lambda_1)^{2+d/2}}\geq 0.
\end{equation}
By the Weyl's law this implies the upper bound
\begin{equation}\label{R-2-Weyl-shifted-upper-bound}
  R_2(z)\leq L_{2,d}^{class}|\Omega|(z+\frac{d}{4}\,\Lambda_1)^{2+d/2}.
\end{equation}
\end{cor}

When $\Omega=\mathcal M=\mathbb{S}^d$, it has been conjectured in \cite{EHI}, on the basis of computations using a computer algebra system, that
\begin{equation}\label{Polynomial-identity-sphere}
  P_N(z)= Q_N(z)
\end{equation}
for all $N\geq 1$ such that $\lambda_{N}<\lambda_{N+1}$. This was first proven in \cite{R}. We prove here this identity for all CROSS:

\begin{thm}\label{main_thm_3}
For CROSS we have
$$
P_N(z)=Q_N(z)
$$
for all $N\geq 1$ such that $\lambda_N<\lambda_{N+1}$.
\end{thm}
The proof of Theorem \ref{main_thm_3}, contained in Section \ref{proof2}, relies in the explicit knowledge of the multiplicity of the eigenvalues of the Laplacian (see also \cite{hel}). In appendix \ref{AppB} we show that for the equilateral and the square flat $2$-dimensional tori, the identity $P_N=Q_N$ fails.

\medskip

\section{Proof of Theorem \ref{main_thm_2}}\label{proof1}

We start with the  algebraic identity for the spectrum of an abstract selfadjoint operator $H$ defined on a Hilbert space $\mathcal{H}$ with scalar product $\langle \cdot,\cdot\rangle_{\mathcal H}$  shown in  \cite{HaSt11}, which we restate here in the particular case that the spectrum of $H$  consists of eigenvalues $\lambda_j$, with an
orthonormal basis of eigenfunctions
$\left\{\phi_j\right\}$, $j\in\mathbb N$.  Let $\mathcal D_H$ be the domain of definition of $H$  and let $G$ be a linear operator with domain $\mathcal D_G$ and adjoint $G^*$ defined on $\mathcal D_{G^*}$. Assume that $G(\mathcal D_H)\subseteq\mathcal D_H\subseteq\mathcal D_G$ and $G^*(\mathcal D_H)\subseteq\mathcal D_H\subseteq\mathcal D_{G^*}$. Let $J$ be a subset of the spectrum of $H$. Then
\begin{equation}\label{HS-sum-rule-abstract}
\begin{split}
    &\;\frac1{2}\sum_{\lambda_j\in J}  (z - \lambda_j)^2\,\big(\langle[G^*,[H,G]]\phi_j,\phi_j\rangle_{\mathcal H}+\langle[G, [H,G^*]]\phi_j,\phi_j\rangle_{\mathcal H}\big)\\
    &-\sum_{\lambda_j\in
    J}(z-\lambda_j)\,\big(\langle[H,G]\phi_j,[H,G]\phi_j\rangle_{\mathcal H}+\langle[H,G^*]\phi_j,[H,G^*]\phi_j\rangle_{\mathcal H}\big)\\
    &=\\
    &\sum_{\lambda_j\in J}\sum_{\lambda_k\notin J}
    (z-\lambda_j)(z-\lambda_k)(\lambda_k-\lambda_j)\big(|\langle G\phi_j,\phi_k\rangle_{\mathcal H}|^2+|\langle G^*\phi_j,\phi_k\rangle|_{\mathcal H}^2\big),\\
    \end{split}
\end{equation}
Here $[H,G]:=HG-GH$ denotes the first commutator of $H$ and $G$, and the first commutator $[H,G^*]$ and the second commutators $[G,[H,G^*]]$ and $[G^*,[H,G]]$ are defined analogously.

If $J=\{\lambda_1,\ldots,\lambda_N\}$ then the right-hand side of \eqref{HS-sum-rule-abstract} has a sign for all $z\in [\lambda_N,\lambda_{N+1}]$. More precisely, an upper bound for \eqref{HS-sum-rule-abstract} on this interval is given by
\begin{equation}\label{HS-sum-rule-abstract-upper-bound}
  \frac1{2}\,(z - \lambda_N)(z - \lambda_{N+1})\sum_{\lambda_j\in J} \,\left(\langle[G^*,[H,G]]\phi_j,\phi_j\rangle_{\mathcal H}+\langle[G, [H,G^*]]\phi_j,\phi_j\rangle_{\mathcal H}\right)
\end{equation}
To prove Theorem \ref{main_thm_2}, we will use \eqref{HS-sum-rule-abstract} (and \eqref{HS-sum-rule-abstract-upper-bound}) with $H=-\Delta$, the Laplacian on $\mathcal M$ or the Dirichlet Laplacian on a domain of $\Omega\subset \mathcal M$,  $\mathcal H=L^2(\mathcal M)$ or $L^2(\Omega)$ with its standard scalar product $\langle f,g\rangle_{\mathcal H}:=\int_{\Omega}fg\,dv$, where $dv$ is the Riemannian volume form. We choose $G$ as a multiplication operator by a suitable function (which we still indicate by $G$). For simplicity we may suppose that $G$ is a real-valued function, that is $G=G^*$. In particular, we take $G$ to be an eigenfunction of the Laplacian on $\mathcal M$: $-\Delta G= \Lambda G$, $\Lambda\geq 0$.  Then
\begin{equation}\label{H-G-commutators-real-case}
  [H,G]=-\Delta G -2\langle\nabla G ,\nabla\cdot\rangle= \Lambda G -2\langle\nabla G,\nabla \cdot\rangle,\quad [G,[H,G]]=2|\nabla G|^2.
\end{equation}
Here by $\langle\cdot,\cdot\rangle$ we denote the scalar product induced by the metric of $\mathcal M$ on tangent spaces. Let now $u\in \mathcal D_H$ (recall that $\mathcal D_H=H^2(\mathcal M)$ or $H^1_0(\Omega)\cap H^2(\Omega)$). Integrating by parts

\begin{equation}\label{H-G-commutator-expectation0}
\int_{\Omega}([G,[H,G]]u)u\,dv=2\int_{\Omega}|\nabla G|^2|u|^2\,dv
\end{equation}
 \begin{equation}\label{H-G-commutator-expectation}
  \int_{\Omega}|[H,G]u|^2\,dv= \int_{\Omega}\left(2\Lambda|\nabla G|^2|u|^2-\Lambda^2|G|^2|u|^2\right)\,dv+4\int_{\Omega}|\langle\nabla G,\nabla u\rangle|^2\,dv.
\end{equation}

Taking $u=\phi_j$ in \eqref{H-G-commutator-expectation0} and \eqref{H-G-commutator-expectation}, we find that the sum rule \eqref{HS-sum-rule-abstract} reads as follows:
\begin{equation}\label{HS-sum-rule-Laplace-Beltrami}
\begin{split}
    &\;\sum_{\lambda_j\in J}  (z - \lambda_j)^2\,\int_{\Omega}|\nabla G|^2|\phi_j|^2\,dv\\
    &-\sum_{\lambda_j\in
    J}(z-\lambda_j)\,\left(\int_{\Omega}\left(2\Lambda|\nabla G|^2|\phi_j|^2-\Lambda^2|G|^2|\phi_j|^2\right)\,dv+4\int_{\Omega}|\nabla G\cdot\nabla \phi_j|^2\,dv\right)\\
    &=\\
    &\sum_{\lambda_j\in J}\sum_{\lambda_k\notin J}
    (z-\lambda_j)(z-\lambda_k)(\lambda_k-\lambda_j)\left(|\langle G\phi_j,\phi_k\rangle_{\mathcal H}|^2\right),\\
    \end{split}
\end{equation}
As already noticed before, if $J=\{\lambda_1,\ldots,\lambda_N\}$ then the right-hand side of the above identity has a sign for all $z\in [\lambda_N,\lambda_{N+1}]$ and  according to \eqref{HS-sum-rule-abstract-upper-bound} an upper bound is given by
\begin{equation}\label{HS-sum-rule-Laplacian-upper-bound}
  (z - \lambda_N)(z - \lambda_{N+1})\sum_{\lambda_j\in J} \,\int_{\Omega}|\nabla G|^2|\phi_j|^2\,dv.
\end{equation}
We choose $G(x)= Y_1^{\alpha}(x)$ for each $\alpha=1,\ldots,M_1$, and hence $\Lambda=\Lambda_1$. Altogether we have
\begin{equation}\label{before_average}
\begin{split}
    &\;\sum_{j=1}^N  (z - \lambda_j)^2\,\int_{\Omega}|\nabla Y_1^{\alpha}|^2|\phi_j|^2\,dv\\
    &-\sum_{j=1}^N(z-\lambda_j)\,\left(\int_{\Omega}\left(2\Lambda_1|\nabla Y_1^{\alpha}|^2|\phi_j|^2-\Lambda_1^2|Y_1^{\alpha}|^2|\phi_j|^2\right)\,dv+4\int_{\Omega}|\nabla Y_1^{\alpha}\cdot\nabla \phi_j|^2\,dv\right)\\
&=\sum_{j=1}^N\sum_{k=N+1}^{\infty}
    (z-\lambda_j)(z-\lambda_k)(\lambda_k-\lambda_j)\left(|\langle Y_1^{\alpha}\phi_j,\phi_k\rangle_{\mathcal H}|^2\right)\\
    &\leq 
    (z - \lambda_N)(z - \lambda_{N+1})\sum_{j=1}^N \,\int_{\Omega}|\nabla Y_1^{\alpha}|^2|\phi_j|^2\,dv.
    \end{split}
\end{equation}
Finally, we average \eqref{before_average} over $\alpha=1,...,M_1$. Then, using identities \eqref{addition-formula2} and \eqref{addition-formula3} for Laplacian eigenfunctions in compact homogeneous manifolds , and identity \eqref{addition-formula5} from Theorem \ref{main_thm}, we prove Theorem \ref{main_thm_2}.

\begin{rem}
Theorem \ref{main_thm_3} states that, when $\mathcal M$ is a CROSS, then $P_N(z)=Q_N(z)$ for all $N\geq 1$ such that $\lambda_N<\lambda_{N+1}$. Inspecting the proof of Theorem \ref{main_thm_2} above, we find that in the case $\Omega=\mathcal M$ is a CROSS:
\begin{multline}\label{selection_rule}
\sum_{j=0}^{\ell}\sum_{k=\ell+1}^{\infty}\sum_{m=1}^{M_j}\sum_{n=1}^{M_k}(z-\Lambda_j)(z-\Lambda_k)(\Lambda_k-\Lambda_j)\sum_{\alpha=1}^{M_1}\left(\int_{\mathcal M}Y_1^{\alpha}Y_j^mY_k^n\,dv\right)^2\\
=(z-\Lambda_{\ell})(z-\Lambda_{\ell+1})\frac{M_1\Lambda_1}{|\mathcal M|}\sum_{j=0}^{\ell}M_j.
\end{multline}
More precisely, averaging \eqref{before_average} over $\alpha=1,...,M_1$, we get $P_N(z)\leq Q_N(z)$. In the case that $\Omega=\mathcal M$ is a CROSS we have equality, which means that the (averaged) inequality between the third and the fourth lines of \eqref{before_average} is an equality. Exploiting the structure of the spectrum of the whole manifold $\mathcal M$, that is, a set of numbers $\Lambda_{\ell}$ with multiplicities $M_{\ell}$, we get \eqref{selection_rule}. Since \eqref{selection_rule} holds for all $\ell$, we deduce that
\begin{equation}\label{selection_rule_2}\int_{\mathcal M}Y_1^{\alpha}Y_j^mY_k^n\, dv=0
\end{equation}
for all $j,k$ such that $|k-j|>1$. This identity is a generalization of what is known in quantum mechanics as a {\it selection rule}. Indeed, in the case of spherical harmonics it corresponds to the transition dipole moment for a single particle state changing from $Y_j^m$ to $Y_k^n$ (see e.g.  \cite[Chapter IV]{bethe_book} where $Y_1^{\alpha}$ represents the angular part of the coordinate function $X_{\alpha}$). Such integrals are special cases of the so-called Slater integrals (see \cite{slater_book}).
\end{rem}

\section{Proof of Theorem \ref{main_thm_3}}\label{proof2}

In this section we prove \eqref{Polynomial-identity-sphere} for all CROSS. To do so,  we first establish conditions yielding a basic sum rule for sequences of numbers, then we apply these results to the spectrum of the CROSS.

\subsection{A basic sum rule for sequences}\label{recurrence_proof}
 Let $\{\lambda_j\}_{j=1}^{\infty}$ be a non-decreasing sequence of non-negative numbers. Let $a>1$ and $h>0$ be positive numbers that will be specified later. For any positive integer $N$ we define the following quadratic polynomials
\begin{equation}\label{P-N}
    P_N(z):=\sum_{j=1}^{N}(z-\lambda_j)(z-h-a\lambda_j)
\end{equation}
and
\begin{equation}\label{Q-N}
    Q_N(z):=N(z-\lambda_N)(z-\lambda_{N+1}).
\end{equation}
Consequently,
\begin{equation}\label{P-N-diff}
    P_{N+1}(z)-P_N(z)=(z-\lambda_{N+1})(z-h-a\lambda_{N+1})
\end{equation}
and
\begin{equation}\label{Q-N-diff}
   Q_{N+1}(z)-Q_N(z)=(z-\lambda_{N+1})(z-\lambda_{N+2}-N(\lambda_{N+2}-\lambda_{N+1}))
\end{equation}
We want to establish conditions ensuring that $\displaystyle P_N(\lambda_{N})=P_N(\lambda_{N+1})=0$ for all positive integers $N$. In particular, if there is a gap between
$\lambda_{N}$ and $\lambda_{N+1}$, that is, if $\lambda_{N}<\lambda_{N+1}$, then $P_N(z)=Q_N(z)$.
Since $P_N(\lambda_{N+1})=0$ implies $P_{N+1}(\lambda_{N+1})=0$ by \eqref{P-N-diff}, we will proceed by induction.

\medskip

First of all, we note that without loss of generality we may assume $h=0$ and $\lambda_1>0$. In fact, let $\tilde\lambda_j:=\lambda_j+\frac{h}{a+1}>0$ and $\tilde z:=z-\frac{h}{a-1}$. Then
$$
\sum_{j=1}^N(z-\tilde\lambda_j)(z-a\tilde\lambda_j)=\sum_{j=1}^N(\tilde z-\lambda_j)(\tilde z-h-a\lambda_j).
$$

Hence it is enough to consider the case $h=0$ and $\lambda_1>0$. We have that
\begin{equation}\label{P-1}
    P_1(z)=(z-\lambda_1)(z-a\lambda_1)
\end{equation}
and
\begin{equation}\label{Q-1}
    Q_1(z)=(z-\lambda_1)(z-\lambda_{2}).
\end{equation}
Then $P_1(\lambda_1)=0$. If $\lambda_2>\lambda_1$ we choose $a=\lambda_2/\lambda_1>1$ in order to get $P_1(\lambda_2)=0$.  If $\lambda_2=\lambda_1$ obviously $P_1(\lambda_2)=0$ and we proceed until the first positive integer $m_0$ such that $\lambda_{m_0+1}>\lambda_{m_0}$ and choose $a=\lambda_{m_0+1}/\lambda_{m_0}>1$. Then $P_{m_0}(\lambda_{m_0+1})=P_{m_0}(\lambda_{m_0+1})=0$ and consequently $P_{m_0}(z)=Q_{m_0}(z)$. Therefore we also may assume without loss of generality $\lambda_2>\lambda_1$ and $a=\lambda_2/\lambda_1>1$.
%\begin{rem}
%  For $h\neq 0$ this conditions transforms to $h=\lambda_2-a\lambda_1$.
%\end{rem}

\medskip

Now let  $\lambda_{N}<\lambda_{N+1}$. By induction hypothesis We may suppose $P_{N}(\lambda_{N})=0$. We compute
\begin{equation*}
 P_{N}(\lambda_{N+1})=P_{N}(\lambda_{N+1})-P_{N}(\lambda_{N})=(\lambda_{N+1}-\lambda_{N})\sum_{j=1}^{N}(\lambda_{N+1}+\lambda_{N}-(a+1)\lambda_j)
\end{equation*}
We therefore have the condition
\begin{equation}\label{Condition-general}
  N(\lambda_{N+1}+\lambda_{N})=(a+1)\sum_{j=1}^{N}\lambda_j
\end{equation}
whenever $\lambda_{N}<\lambda_{N+1}$.

\subsection{Basic example: strictly increasing sequences} We suppose $0<\lambda_{N}<\lambda_{N+1}$ for all $N$. Since by definition $a=\frac{\lambda_2}{\lambda_1}$, condition \eqref{Condition-general} holds for $N=1$.
Suppose now that condition \eqref{Condition-general} holds for some $N$. Then by hypothesis
\begin{equation*}
\begin{split}
  & (N+1)(\lambda_{N+2}+\lambda_{N+1})-(a+1)\sum_{j=1}^{N+1}\lambda_j\\
  &=(N+1)(\lambda_{N+2}+\lambda_{N+1})-N(\lambda_{N+1}+\lambda_{N})-(a+1)\lambda_{N+1}\\
  &=(N+1)\lambda_{N+2}-a\lambda_{N+1}-N\lambda_{N}.\\
\end{split}
\end{equation*}
A solution of this recurrence relation can be found for particular values of $a$. For example, if $a=3$ then $\lambda_{N}=(2N-1)\lambda_{1}$ satisfies the recurrence relation; if $a=5$ then
$\lambda_{N}=(2N^2-2N+1)\lambda_{1}$ satisfies the recurrence relation.
\subsection{General case: sequences with multiplicities} Let $m_j\geq 0$ be positive integers and $\{\Lambda_j\}_{j=0}^{\infty}$ be a strictly increasing sequence of non-negative numbers. We define the sequence $\lambda_j$ by
\begin{equation*}
  \lambda_1=\cdots=\lambda_{m_0}=\Lambda_0\,,\ \ \  \lambda_{m_0+1}=\cdots = \lambda_{m_0+m_1}=\Lambda_1\,,\ \ \ \cdots
\end{equation*}
and the number $\#\{\lambda_j:\lambda_j\leq \Lambda_n\}$ by 
\begin{equation}\label{N-n-def}
  N_n:=\sum_{j=0}^{n}m_j.
\end{equation}
Following Subsection \ref{recurrence_proof} in the case $\lambda_{N_n}<\lambda_{N_n+1}$, which is $\Lambda_{n}<\Lambda_{n+1}$, we find the following condition for the sequence
\begin{equation}\label{Lambda_tilde}
  \tilde{\Lambda}_n=\Lambda_n+\frac{h}{a-1}:
\end{equation}
\begin{equation}\label{Condition-energy-values}
 N_{n+1}( \tilde{\Lambda}_{n+2}-a \tilde{\Lambda}_{n+1})+ N_{n}( a\tilde{\Lambda}_{n+1}- \tilde{\Lambda}_{n})=0,
\end{equation}
where
$$
a=\frac{\tilde\Lambda_1}{\tilde\Lambda_0}
$$
Observe that in order to have solutions of this recurrence relation all $\tilde{\Lambda}_{n}$ must satisfy the growth condition
\begin{equation}\label{Growth-Condition-energy-values}
 \tilde{\Lambda}_{n+2}-a \tilde{\Lambda}_{n+1}<0.
\end{equation}
We can rewrite the relation \eqref{Condition-energy-values} as a first order recurrence relation for $N_n$ which we can solve explicitly. However, we always have to check that the solutions $N_n$ are positive integers  strictly increasing in $n$. Namely, we want to solve the equation
\begin{equation}\label{N-recurrence}
  N_{n+1}= \frac{a\tilde{\Lambda}_{n+1}- \tilde{\Lambda}_{n}}{a\tilde{\Lambda}_{n+1}- \tilde{\Lambda}_{n+2}}\,N_n.
\end{equation}
We see indeed that always $N_{n+1}>N_{n}$.

\subsection{Harmonic oscillator: $\mathbf{\Lambda_l=l+1/(a-1)}$} With this choice $h=0$ and  $\tilde{\Lambda}_{l}=\Lambda_{l}$. The recurrence relation \eqref{N-recurrence} becomes
\begin{equation*}
  N_{l+1}=\left(1+\frac{2}{(a-1)(l+1)}\right)\,N_l
\end{equation*}
which admits the solution
\begin{equation*}
  N_l=\frac{\Gamma(l+1+\frac{2}{a-1})}{\Gamma(1+\frac{2}{a-1})\Gamma(l+1)}\,N_0=\binom{l+\frac{2}{a-1}}{l}\,N_0.
\end{equation*}
For $a=2$ and $a=3$ and $N_0=1$, the sequences $\Lambda_l$ with multiplicities $N_l$ correspond to sequences of eigenvalues (shifted by $-\frac{1}{2}$) of quantum harmonic oscillators in dimension $3$ and $2$, respectively.

\subsection{CROSS: $\mathbf{\Lambda_l=l(l+h-1)}$} The eigenvalues $\Lambda_l$ of CROSS and their multiplicities $m_l$ are know, see \cite{hel}. In this case $\tilde{\Lambda}_l=\Lambda_l+\frac{h}{a-1}$. The recurrence relation \eqref{N-recurrence} becomes
\begin{equation*}
  N_{l+1}=\left(1+\frac{2(2l+1+h)}{(l+1)((a-1)(l+h)-2)}\right)\,N_l
\end{equation*}
which admits the solution when $(a-1)h>2$:
\begin{equation*}
  N_l=\frac{\Gamma(h-\frac{2}{a-1})\Gamma(l+1+\frac{2}{a-1})\Gamma(l+h)}{\Gamma(h)\Gamma(l+\frac{2}{a-1})\Gamma(l+h-\frac{2}{a-1})\Gamma(l+1)}\,N_0
\end{equation*}
or using binomial coefficients
\begin{equation*}
  N_l=\frac{\binom{l+h-1}{l}\binom{l+\frac{2}{a-1}}{l}}{\binom{l+h-1-\frac{2}{a-1}}{l}}\,N_0.
\end{equation*}
For specific choices of $h,a$ we have the spectrum of CROSS:
\begin{itemize}
\item{\bf $d$-dimensional sphere.} Here $h=d$, $a=1+4/d$, $N_0=1$. Hence after simplification for $l\geq 0$:
\begin{equation*}
  N_l=(2+\frac{d}{l})\binom{d+l-1}{d}
\end{equation*}
or for all $l$
\begin{equation*}
  N_l=\frac{(d+2l) \Gamma(l+d)}{\Gamma(d+1)\Gamma(l+1)},
\end{equation*}
which yields the well-known expression
$$
m_l=\frac{2l+d-1}{l}\binom{l-2+d}{l-1}.
$$
\item {\bf $d$-dimensional real projective space.} Here $h=(d+1)/2$, $a=1+4/d$, $N_0=1$. Hence after simplification for $l\geq 0$:
\begin{equation*}
  N_l=\binom{d+2l}{d}
\end{equation*}
which yields for the multiplicities the well-known expression
\begin{equation*}
  m_l=(2+\frac{d-1}{2l})\binom{d+2l-2}{d-1} =\frac{(d-1+4l) \Gamma(2l+d-1)}{\Gamma(d)\Gamma(2l+1)}.
\end{equation*}

\item {\bf $d$-dimensional complex projective space.} Here $h=1+d/2$, $a=1+4/d$, $N_0=1$. Hence after simplification for $l\geq 0$:
\begin{equation*}
  N_l=\binom{l+d/2}{l}^2
\end{equation*}
which yields for the multiplicities the well-known expression
\begin{equation*}
  m_l=(1+\frac{4l}{d})\binom{l-1+d/2}{l}^2.
\end{equation*}

\item {\bf $d$-dimensional quaternion projective space.} Here $h=\frac{d}{2}+2$, $a=1+4/d$, $N_0=1$. For $l\geq 0$ we find
\begin{equation*}
  N_l=\frac{(d+2l+2)(d+2l)}{2dl(l+1)}\binom{l-1+d/2}{l}\binom{l+d/2}{l-1}
\end{equation*}
which yields for the multiplicities the well-known expression
\begin{equation*}
  m_l=\frac{d+4l+2}{2l(l+1)}\binom{l-1+d/2}{l}\binom{l+d/2}{l-1}
\end{equation*}

\item {\bf $d$-dimensional Cayley projective space.} Here $h=\frac{d}{2}+4$, $a=1+4/d$, $N_0=1$. For $l\geq 0$ we find
\begin{equation*}
  N_l=\frac{3(d+2l+6)(d+2l)}{dl(l+1)(l+2)(l+3)}\binom{l-1+d/2}{l}\binom{l+d/2+2}{l-1}
\end{equation*}
which yields for the multiplicities the well-known expression
\begin{equation*}
  m_l=\frac{3(d+4l+6)}{l(l+1)(l+2)(l+3)}\binom{l-1+d/2}{l}\binom{l+d/2+2}{l-1}
\end{equation*}
\end{itemize}

\appendix

\section{Proof of Theorem \ref{main_thm}}\label{AppA}
We provide here a detailed proof of Theorem \ref{main_thm}. We  need to recall a few basic information on compact homogeneous manifolds with irreducible linear isotropy group.

\subsection{Compact homogeneous manifolds with irreducible linear isotropy group}

Let $\mathcal M$ be a compact homogeneous Riemannian manifold of dimension $d\geq 2$. That is, there is a compact Lie group $G$ acting transitively and isometrically on $\mathcal M$, and the metric is $G$-invariant. We denote by $\phi$ the group action. Namely $\phi:G\times \mathcal M\to\mathcal M$  is given by $\phi(g,x)=\phi_g(x)\in \mathcal M$, where the map $x\mapsto \phi_g(x)$ is an isometry. With abuse of notation we will write $g$ in place of $\phi_g$, namely $\phi_g(x)=:g(x)$.

Let $x\in \mathcal M$ and let $G_x$ the stabilizer of $x$, namely
$$
G_x:=\{g\in G:\phi(g,x)=\phi_g(x)=x\}.
$$
We have that
$$
dg_x:T_x\mathcal M\to T_x\mathcal M
$$
is an isometry, which means that
$$
\langle U,V\rangle=\langle dg_xU,dg_xV\rangle
$$
for all $U,V\in T_x\mathcal M$. Through all this section $\langle\cdot,\cdot\rangle$ shall denote the scalar product on $T_x\mathcal M$ associated with the ($G$-invariant) Riemannian metric on $\mathcal M$.

Let $g_x:\mathcal M\to \mathcal M$ be an isometry fixing $x\in \mathcal M$. Consider the map $\rho: G_x\to O(T_x\mathcal M,\langle\cdot,\cdot\rangle)$ given by 
$$
\rho(g_x):=dg_x.
$$
This map is the so-called {\it isotropy representation} at $x$. We give now the definition of {\it irreducible representation}.
\begin{defn}
A subspace $F$ of $T_x\mathcal M$ is called $G$-invariant if $\rho(g_x)V\in F$ for all $g_x\in G_x$ and $V\in F$. We say that $\rho: G_x\to O(T_x\mathcal M,\langle\cdot,\cdot\rangle)$ is irreducible if the only $G$-invariant subspaces of $T_x\mathcal M$ for $\rho$ are $T_x\mathcal M$ and $\{0\}$.
\end{defn}
We can now define an {\it isotropy irreducible Riemannian manifold}.
\begin{defn}\label{def_iso_irr}
A compact homogeneous Riemannian manifold is said {\it isotropy} {\it irreducible} if the isotropy representation is irreducible on tangent spaces.
\end{defn}
It is well-known that all CROSS are isotropy irreducible Riemannian manifolds (they are isotropic). We refer to \cite[\S I.F,\S III.C]{BGM}, \cite{wolf} and \cite[\S 8.13]{wolf_book} for more details and examples on isotropy irreducible homogeneous spaces. Note that an isotropy irreducible Riemannian manifold can be written as a {\it homogeneous space} $G/G_x$, see e.g., \cite[\S 1]{wang_ziller}.

\subsection{Proof of Theorem \ref{main_thm}}

We are now in position to prove Theorem \ref{main_thm}. Let $\Lambda_l>0$ be an eigenvalue level, and let $Y_l^k$, $k=1,...,M_l$ be a $L^2(\mathcal M)$ orthonormal basis of real-valued eigenfunctions of the corresponding eigenspace. Since the Laplacian commutes with the isometries we have that for any eigenfunction $Y_k^l$ associated with $\Lambda_l$, then $Y_l^k\circ g$ is another eigenfunction associated with the same eigenvalue.

Consider the linear operator $L:T_x\mathcal M\to T_x\mathcal M$ defined by
$$
LV=\sum_{k=1}^{M_l}\langle V,\nabla Y_l^k(x)\rangle\nabla Y_l^k(x).
$$
From now on we suppress the dependence of $\nabla Y_l^k(x)$ from $x$ when it is clear from the context. We will prove that $L$ is an homothety, namely, we will prove that
\begin{equation}\label{hom}
LV=\frac{M_l\Lambda_l}{d|\mathcal M|}V
\end{equation}
for all $V\in T_x\mathcal M$.
 Let  $g_x\in G_x$. We compute
\begin{multline}\label{part1}
L(dg_xV)=\sum_{k=1}^{M_l}\langle dg_xV,\nabla Y_l^k\rangle\nabla Y_l^k=\sum_{k=1}^{M_l}\langle V, d_{g_x^{-1}}\nabla Y_l^k\rangle\nabla Y_l^k\\
=dg_x\sum_{k=1}^{M_l}\langle V,\nabla (Y_l^k\circ g_x^{-1})\rangle\nabla (Y_l^k\circ g_x^{-1}).
\end{multline}
Here we have used the fact that $dg_x\in O(T_x\mathcal M,\langle\cdot,\cdot\rangle)$ and that $\nabla(Y_l^k\circ g)=dg\nabla Y_l^k$ for any $g\in G_x$. Now, since $g_x^{-1}\in G_x$, we have that $Y_k^l\circ g_x^{-1}$ is an eigenfunction associated with $\Lambda_l$, since the eigenspace is invariant under the action of $G$. Therefore, setting $\tilde Y_l^k:=Y_l^k\circ g_x^{-1}$, we have
$$
\tilde Y_l^k=\sum_{j=1}^{M_l}a_{kj}Y_l^j.
$$
Since the metric is invariant under isometries, we have that $\int_{\mathcal M}\tilde Y_l^k\tilde Y_l^j=\delta_{kj}$. This implies that
$$
\delta_{kj}=\int\tilde Y_l^j\tilde Y_l^j=\sum_{m,n}\int_M a_{km}a_{j n}Y_{l}^mY_l^n=\sum_{m=1}^{M_l}a_{km}a_{jm},
$$
which means that the matrix $(a_{kj})$ is orthogonal. Going back to \eqref{part1} we have
\begin{equation}\label{part2}
L(dg_xV)=dg_x\sum_{k=1}^{M_l}\sum_{j,m=1}^{M_l}a_{kj}a_{km}\langle V,\nabla Y_k^l\rangle\nabla Y_m^l=dg_x\sum_{k=1}^{M_l}\langle V, Y_l^k\rangle\nabla Y_l^k=dg_xL(V).
\end{equation}
Therefore, $L$ is a linear operator on $T_x\mathcal M$ commuting with all $g_x\in G_x$ (or, to be more precise, with the image of $G_x$ in $O(T_x\mathcal M,\langle\cdot,\cdot\rangle)$ by the representation $\rho$).

Moreover, by definition, $L$ is symmetric:
$$
\langle L(U),V\rangle=\langle U,L(V)\rangle
$$
for all $U,V\in T_x\mathcal M$. Hence $T_x\mathcal M$ admits an orthonormal basis of eigenfunctions of $L$. Let $F_{\mu}$ be an eigenspace of $L$ associated with an eigenvalue $\mu$. Consider the space $\{dg_xV:g_x\in G_x,V\in F_{\mu}\}$. Since for $V\in F_{\mu}$ we have $L(dg_xV)=dg_xL(V)=\mu dg_xV$, we conclude that $dg_xV$ is an eigenfunction associated with $\mu$. Hence $F_{\mu}$ is a $G$-invariant subspace of $T_x\mathcal M$, and from the hypothesis of irreducibility, necessarily $F_{\mu}=T_x\mathcal M$ (since it is not the trivial subspace $\{0\}$). This implies that there is a unique eigenvalue $\mu\in\mathbb R$ and $L=\mu\, Id$ on $T_x\mathcal M$.

The previous computations also show that $\mu$ does not depend on $x$. \\
In fact, taking another point $x'$, we can define $L':T_{x'}\mathcal M\to T_{x'}\mathcal M$ by $L'=\sum_{i=k}^{M_l}\langle V,\nabla Y_l^k(x')\rangle\nabla Y_l^k(x')$. Here and in what follows the scalar product is the one on the appropriate tangent space. Now, there exists $g\in G$ such that $x'=g(x)$. Then
\begin{multline*}
L'(V)=\sum_{k=1}^{M_l}\langle V,\nabla Y_l^k(g(x))\rangle\nabla Y_l^k(g(x))=dg\sum_{k=1}^{M_l}\langle V, dg\nabla(Y_l^k\circ g)\rangle\nabla(Y_l^k\circ g)\\
=dg\sum_{k=1}^{M_l}\langle dg^{-1}V,\nabla(Y_l^k\circ g)\rangle\nabla(Y_l^k\circ g)=dg\sum_{k=1}^{M_l}\langle dg^{-1}V,\nabla Y_l^k\rangle\nabla Y_l^k\\
=dgL(dg^{-1}V)=\mu V.
\end{multline*}

To conclude the proof, we need to compute the value of $\mu$. To do so, let $(e_1,...,e_d)$ be a orthonormal basis of $T_x\mathcal M$. We have
$$
\langle L(e_i),e_i\rangle=\mu
$$
which reads
$$
\sum_{k=1}^{M_l}|\langle e_i,\nabla Y_k^l\rangle|^2=\mu.
$$
Summing over $i$ we get
$$
\sum_{k=1}^{M_l}|\nabla Y_l^k|^2=\sum_{k=1}^{M_l}\sum_{i=1}^d|\langle e_i, Y_l^k\rangle|^2=\mu d.
$$
In particular $\mu\ne 0$ (none of the $Y_k^l$ is constant). Integrating over $\mathcal M$ we get
$$
M_l\Lambda_l=\sum_{k=1}^{M_l}\int_{\mathcal M}|\nabla Y_l^k|^2=\int_{\mathcal M}\mu d=|\mathcal M|\mu d.
$$
We conclude that
$$
\mu=\frac{M_l\Lambda_l}{d|\mathcal M|}.
$$
This concludes the proof of \eqref{hom} and of Theorem \ref{main_thm}. \qed

\section{The flat torus}\label{AppB}

A $2$-dimensional flat torus $T$ is defined as $T=\mathbb R^2/\Gamma$ with the induced metric from $\mathbb R^2$. Here $\Gamma$ is a lattice, namely, $\Gamma=\{nw_1+mw_2:n,m\in\mathbb Z\}$ and $(w_1,w_2)$ is a basis of $\mathbb R^2$. The dual lattice $\Gamma^*$ is defined as the set $\Gamma^*:=\{p\in\mathbb R^2:\langle x,p\rangle\in\mathbb Z\ \ \forall x\in\Gamma\}$. 

The eigenfunctions of $T$ are given by $f(x):=e^{2\pi i\langle p,x\rangle}$ where $p$ ranges in $\Gamma^*$. The corresponding eigenvalues are $\lambda=4\pi^2|p|^2$. The multiplicity of a non-zero eigenvalue is always even: if $e^{2\pi i\langle p,x\rangle}$ is an eigenfunction, then also $e^{-2\pi i\langle p,x\rangle}$ is an eigenfunction associated with the same eigenvalue. We always have $\lambda_1=0$ ($p=0$), while $\lambda_2>0$. We refer e.g., to \cite[\S III.B]{BGM} for the computation of the spectrum of flat tori.

It is well-known that, up to a homotheties, congruences in ${\bf O}(2)$ and change of basis in ${\bf SL}(2,{\mathbb Z})$, the lattice $\Gamma$ admits a basis $w_1=(1,0)$, $w_2=(a,b)$ with $(a,b)\in\tau:=\{0\leq a\leq 1/2$, $a^2+b^2\geq 1$, $b>0\}$. The region $\tau$ is the so-called moduli space of flat tori, see \cite[\S III.B]{BGM}. A dual basis is given by $w_1^*=(1,-a/b)$, $w_2^*=(0,1/b)$, and the second eigenvalue is always given by $\lambda_2=\frac{4\pi^2}{b^2}$. Note that with this description of flat tori we are fixing the area to be $b$.

A flat torus is a homogeneous Riemannian manifold, but it is not isotropic. Isotropy irreducible flat tori are just the square and the equilateral torus, and one can easily check by hand that Theorem \ref{main_thm} holds for such tori, while in all the remaining cases Theorem \ref{main_thm} fails. These are known facts since the only flat tori admitting isometric minimal embeddings into spheres are the square and the equilateral ones (see e.g., \cite{EI_00}).

In the following we will explore the validity of $P_N(z)\leq Q_N(z)$ for $z\in[\lambda_N,\lambda_{N+1}]$ when there is a gap and also the failure of the identity $P_N(z)=Q_N(z)$ even for the irreducible cases. This suggests that this property is very specific of CROSS.

\subsection{Square torus} This corresponds to $(a,b)=(0,1)$, then $\Gamma^*=\mathbb Z^2$. We have $\lambda_2=4\pi^2$ with multiplicity $4$ and a basis $\{f_i\}_{i=1}^4$ of the eigenspace is given by $\{e^{\pm 2\pi i \langle p_i,x\rangle}\}_{i=1}^4$ with $p_1=e_1$, $p_2=e_2$, $p_3=-e_1$, $p_4=-e_2$, $(e_1,e_2)$ being the canonical basis of $\mathbb R^2$. If $v=(v_1,v_2)$, then $\sum_i\langle v,\nabla \bar f_i\rangle\nabla f_i=8\pi^2 v$, namely \eqref{addition-formula4}. Note that the crucial point is that $\{p_i\}_{i=1}^4$ forms a {\it tight frame} for $\mathbb R^2$.

\begin{rem} One can write the square torus $T$ also as $T=\mathbb S^1_1\times\mathbb S^1_1$ with the product metric, where by $\mathbb S^1_L$ we denote the circle of length $L$. It is known that if $\mathcal M$ is isotropy irreducible, then also $\mathcal M\times\cdots\times\mathcal M$ is, see e.g., \cite[\S1]{wang_ziller}. This is not always the case for a product (for example, in the case of rectangular tori). Concerning the square torus, one could also check directly the definition of irreducibility, since the stabilizer of a point is the dihedral group $D_8$.
\end{rem}

Clearly from Theorem \ref{main_thm_2} it follows that $P_N(z)\leq Q_N(z)$ for al $z\in[\lambda_N,\lambda_{N+1}]$. However, the identity $P_N(z)=Q_N(z)$ already fails for $N=5$ (recall that the first $6$ eigenvalues are $0,4\pi^2,4\pi^2,4\pi^2,4\pi^2, 8\pi^2$. Then $P_5(z)-Q_5(z)=24\pi^2(4\pi^2-z)$ and \eqref{Polynomial-identity-sphere} does not hold.

\subsection{Equilateral torus} This corresponds to $(a,b)=(1/2,\sqrt{3}/2)$. We have $\lambda_2=\frac{16\pi^2}{3}$ with multiplicity $6$ and a basis $\{f_i\}_{i=1}^6$ of the eigenspace is given by $\{e^{\pm 2\pi i \langle p_i,x\rangle}\}_{i=1}^6$ with $p_1=(1,1/\sqrt{3})$, $p_2=(0,2/\sqrt{3})$, $p_3=(-1,1/\sqrt{3})$, $p_4=(-1,-1/\sqrt{3})$, $p_5=(0,-2/\sqrt{3})$, $ p_6=(1,-1/\sqrt{3})$. Again, if $v=(v_1,v_2)$, then $\sum_i\langle v,\nabla \bar f_i\rangle\nabla f_i=16\pi^2 v$. Also in this case we note that $\{p_i\}_{i=1}^6$ forms a tight frame for $\mathbb R^2$.

Note that $\lambda_2$ has multiplicity $6$ which is the maximal possible multiplicity for a flat torus, and is attained only by the equilateral torus. Again, the irreducibility could be checked by hand, noting that the stabilizer of a point in this case is the dihedral group $D_{12}$.

In this case, the first $8$ eigenvalues are $0,\frac{16\pi^2}{3},\frac{16\pi^2}{3},\frac{16\pi^2}{3},\frac{16\pi^2}{3},\frac{16\pi^2}{3},\frac{16\pi^2}{3},16\pi^2$ and $P_7(z)-Q_7(z)=16\pi^2\left(\frac{16\pi^2}{3}-z\right)$. Then \eqref{Polynomial-identity-sphere} does not hold.

\subsection{Rectangular torus} Here we have $(a,b)=(0,b)$ with $b>1$. The dual lattice $\Gamma^*$ is generated by $w_1^*=(1,0)$, $w_2^*=(0,1/b)$. We have $\lambda_2=\frac{4\pi^2}{b^2}$, with multiplicity $2$ and a basis of the eigenspace is $\{e^{2\pi i \langle p_i,x,\rangle}\}$, $p_1=(0,1/b)$, $p_2=(0,-1/b)$. Now, $\{p_1,p_2\}$ does not form a tight frame, hence we don't have \eqref{addition-formula4}.
One can check the non irreducibility also from the definition: the stabilizer of a point is the Klein group (or, equivalently, the product $C_2\times C_2$ of two cyclic groups), hence one can identify  one dimensional invariant subspaces of the tangent space at $x$ for the group action.

Of course, the fact that a torus is not irreducible does not imply that $P_N\leq Q_N$ does not hold. Assume that $b\geq 2$. The first four eigenvalues are given by $0,\frac{4\pi^2}{b^2},\frac{4\pi^2}{b^2},\frac{16\pi^2}{b^2}$. We compute $P_3(z)-Q_3(z)=\frac{16\pi^2}{b^4}(b^2z-4\pi^2)\geq 0$ on $[\frac{4\pi^2}{b^2},\frac{16\pi^2}{b^2}]$. If $1<b<2$,  the first four eigenvalues are $0,\frac{4\pi^2}{b^2},\frac{4\pi^2}{b^2},4\pi^2$ and we see that $P_3(z)-Q_3(z)\leq 0$ if and only if $1<b<2\sqrt{2/3}$. If $1<b<2\sqrt{2/3}$ it is not clear whether there exist some $N$ for which $P_N\geq Q_N$. For these values of $b$ we have $P_3\leq Q_3$ and $P_5\leq Q_5$.

\subsection{In general}

Consider now any $(a,b)\in\tau$. If $a^2+b^2>1$, then the multiplicity of $\lambda_2=\frac{4\pi^2}{b^2}$ is  $2$ and an eigenbasis is given by $\{e^{2\pi i\langle p_i,x\rangle}\}_{i=1}^2$ with $p_1=(0,1/b)$, $p_2=(0,-1/b)$ which is not a tight frame. It remains to consider $a^2+b^2=1$, $0\leq a\leq 1/2$. We have already seen the case $a=0$ (square torus) and $a=1/2$ (equilateral torus). In the first case the multiplicity of $\lambda_2$ is $4$, in the second is $6$. We see that for any $0< a<1/2$ the multiplicity is $4$ and an eigenbasis is given by $\{e^{2\pi i\rangle p_i,x\rangle}\}_{i=1}^4$ with $p_1=(1,-a/b)$, $p_2=(-1,a/b)$, $p_3=(0,1/b)$, $p_4=(0,-1/b)$, which is not a tight frame for any $0<a<b$, hence \eqref{addition-formula4} does not hold. 

Concerning the inequality $P_N\leq Q_N$, we see that for large enough values of $b$, we have that the bottom of the spectrum is given by $0,\frac{4\pi^2}{b^2},\frac{4\pi^2}{b^2},\frac{16\pi^2}{b^2},...$. As in the case of the rectangular torus we see easily that $P_3(z)\geq Q_3(z)$. On the other hand, for small values of $N$ it is possible to check that $P_N\leq Q_N$ whenever $a^2+b^2=1$. In view of this, it is natural to conjecture that for $(a,b)$ in $\{a^2+b^2=1\}\cap\{0\leq a\leq 1/2\}\cap\{b>0\}\times[0,+\infty]$, and in a sufficiently small neighborhood of this set, the inequality $P_N\leq Q_N$ holds for any $N$ when there is a gap.

\subsection*{Alternative definitions of $P_N$}

When treating non irreducible tori, it could be convenient to redefine $P_N$. Consider \eqref{HS-sum-rule-Laplace-Beltrami} and \eqref{HS-sum-rule-Laplacian-upper-bound} and take $G=e^{2\pi i\langle p,x\rangle}$ with $p$ in the dual lattice $\Gamma^*$. Hence $G$ is (the multiplication by) an eigenfunction of the flat torus $\mathbb R^2/\Gamma$ with associated eigenvalue $4\pi^2|p|^2$. Hence we get the inequality
\begin{equation}\label{HS_torus}
\sum_{j=1}^N(z-\lambda_j)\left(z-4\pi^2|p|^2-4\int_{\Omega}\left|\langle\frac{p}{|p|},\nabla\phi_j\rangle\right|^2-\lambda_j\right)\leq N(z-\lambda_N)(z-\lambda_{N+1}).
\end{equation}
We can consider now a basis $p_1,p_2$ of $\Gamma^*$, but more in general we can consider $p_1,...,p_k\in\Gamma^*$. Call $L_i:=4\pi^2|p_i|^2$. Note that these are not the ordered eigenvalues, just some eigenvalues that we label $L_1,...,L_k$. Then averaging \eqref{HS_torus} we get
\begin{multline}\label{HS_torus_2}
\sum_{j=1}^N(z-\lambda_j)\left(z-\frac{1}{k}\sum_{i=1}^kL_i-\frac{4}{k}\int_{\Omega}\sum_{i=1}^k\left|\langle\frac{p}{|p|},\nabla\phi_j\rangle\right|^2-\lambda_j\right)\\
\leq N(z-\lambda_N)(z-\lambda_{N+1}).
\end{multline}

The question is if we find a tight frame which consists of $\frac{p_i}{|p_i|}$. This is in fact the case of rectangular tori, though we don't necessarily have the first two eigenvalues. In fact in this case we can take $p_1=(0,1/b)$, $p_2=(1,0)$, and hence $L_1=\frac{4\pi^2}{b^2}=\Lambda_1$ (the first positive eigenvalue) and $L_2=4\pi^2$ (some eigenvalue, not necessarily the second eigenvalue level). Note that as $b\to\infty$, $L_2=\Lambda_k$ with $k\to\infty$. Inequality \eqref{HS_torus_2} reads for such rectangular tori
\begin{equation}\label{HS_rect_torus}
\sum_{j=1}^N(z-\lambda_j)\left(z-2\pi^2(1+1/b^2)-3\lambda_j\right)\\
\leq N(z-\lambda_N)(z-\lambda_{N+1}).
\end{equation}
Note that one would like to replace the shift $\frac{L_1+L_2}{2}$ by $\frac{\Lambda_1+\Lambda_2}{2}$ which seems more natural. However, taking $b\to+\infty$, we see that $\Lambda_1=\frac{4\pi^2}{b^2}$, $\Lambda_2=\frac{16 p^2}{b^2}$, $\Lambda_3=\frac{36\pi^2}{b^2}$ each of multiplicity  $2$, and $P_5(36\pi^2/b^2)-Q_5(36\pi^2/b^2)=952\pi^4/b^4$.

In general, if $a\in\mathbb Q$, we always find $p_1,p_2\in\Gamma^*$ such that $\langle p_1,p_2\rangle=0$, and hence we have

\begin{equation}\label{HS_torus_3}
\sum_{j=1}^N(z-\lambda_j)\left(z-\frac{L_1+L_2}{2}-3\lambda_j\right) \leq N(z-\lambda_N)(z-\lambda_{N+1}).
\end{equation}
\section*{Acknowledgments}
The first author is grateful to Marco Radeschi for useful discussions on the topic.  The first  author acknowledges support of the project ``Perturbation problems and asymptotics for elliptic differential equations: variational and potential theoretic methods'' funded by the European Union – Next Generation EU and by MUR-PRIN-2022SENJZ. The first author is member of the Gruppo Nazionale per le Strutture Algebriche, Geometriche e le loro Applicazioni (GNSAGA) of the Istituto Nazionale di Alta Matematica (INdAM).

\bibliography{bibliography.bib}

\def\cprime{$'$} \def\cprime{$'$} \def\cprime{$'$} \def\cprime{$'$}
  \def\cprime{$'$}
\begin{thebibliography}{10}

\bibitem{BGM}
M.~Berger, P.~Gauduchon, and E.~Mazet.
\newblock {\em Le spectre d'une vari\'et\'e{} riemannienne}, volume Vol. 194 of
  {\em Lecture Notes in Mathematics}.
\newblock Springer-Verlag, Berlin-New York, 1971.

\bibitem{BJ}
H.~A. Bethe.
\newblock {\em Intermediate quantum mechanics}.
\newblock W. A. Benjamin, Inc., New York-Amsterdam, 1964.
\newblock Notes by R. W. Jackiw.

\bibitem{bethe_book}
H.~A. Bethe and E.~E. Salpeter.
\newblock {\em Quantum mechanics of one- and two-electron atoms}.
\newblock Springer-Verlag, Berlin-G\"ottingen-Heidelberg; Academic Press, Inc.,
  New York, 1957.

\bibitem{EHI}
A.~El~Soufi, E.~M. Harrell, II, and S.~Ilias.
\newblock Universal inequalities for the eigenvalues of {L}aplace and
  {S}chr\"{o}dinger operators on submanifolds.
\newblock {\em Trans. Amer. Math. Soc.}, 361(5):2337--2350, 2009.

\bibitem{EI_00}
A.~El~Soufi and S.~Ilias.
\newblock Riemannian manifolds admitting isometric immersions by their first
  eigenfunctions.
\newblock {\em Pacific J. Math.}, 195(1):91--99, 2000.

\bibitem{G75}
E.~Gin\'{e}~M.
\newblock The addition formula for the eigenfunctions of the {L}aplacian.
\newblock {\em Advances in Math.}, 18(1):102--107, 1975.

\bibitem{Ha0}
E.~M. Harrell, II.
\newblock Commutators, eigenvalue gaps, and mean curvature in the theory of
  {S}chr\"odinger operators.
\newblock {\em Comm. Partial Differential Equations}, 32(1-3):401--413, 2007.

\bibitem{HaHe0}
E.~M. Harrell, II and L.~Hermi.
\newblock Differential inequalities for {R}iesz means and {W}eyl-type bounds
  for eigenvalues.
\newblock {\em J. Funct. Anal.}, 254(12):3173--3191, 2008.

\bibitem{HaSt0}
E.~M. Harrell, II and J.~Stubbe.
\newblock On trace identities and universal eigenvalue estimates for some
  partial differential operators.
\newblock {\em Trans. Amer. Math. Soc.}, 349(5):1797--1809, 1997.

\bibitem{HaSt11}
E.~M. Harrell, II and J.~Stubbe.
\newblock Trace identities for commutators, with applications to the
  distribution of eigenvalues.
\newblock {\em Trans. Amer. Math. Soc.}, 363(12):6385--6405, 2011.

\bibitem{hel}
S.~u. Helgason.
\newblock The {R}adon transform on {E}uclidean spaces, compact two-point
  homogeneous spaces and {G}rassmann manifolds.
\newblock {\em Acta Math.}, 113:153--180, 1965.

\bibitem{Li81}
P.~Li.
\newblock Minimal immersions of compact irreducible homogeneous {R}iemannian
  manifolds.
\newblock {\em J. Differential Geometry}, 16(1):105--115, 1981.

\bibitem{R}
M.~{R}uano.
\newblock {L}es valeurs propres du {L}aplacien sur la sph\`{e}re.
\newblock {\em Master project, EPFL}, 2012.

\bibitem{slater_book}
J.~Slater.
\newblock {\em Quantum Theory of Atomic Structure. Vol. 1}.
\newblock McGraw-Hill, 1960.

\bibitem{St0}
J.~Stubbe.
\newblock Universal monotonicity of eigenvalue moments and sharp
  {L}ieb-{T}hirring inequalities.
\newblock {\em J. Eur. Math. Soc. (JEMS)}, 12(6):1347--1353, 2010.

\bibitem{takahashi}
T.~Takahashi.
\newblock Minimal immersions of {R}iemannian manifolds.
\newblock {\em J. Math. Soc. Japan}, 18:380--385, 1966.

\bibitem{wang_ziller}
M.~Wang and W.~Ziller.
\newblock On isotropy irreducible {R}iemannian manifolds.
\newblock {\em Acta Math.}, 166(3-4):223--261, 1991.

\bibitem{wolf}
J.~A. Wolf.
\newblock The goemetry and structure of isotropy irreducible homogeneous
  spaces.
\newblock {\em Acta Math.}, 120:59--148, 1968.

\bibitem{wolf_book}
J.~A. Wolf.
\newblock {\em Spaces of constant curvature}.
\newblock AMS Chelsea Publishing, Providence, RI, sixth edition, 2011.

\end{thebibliography}
\bibliographystyle{abbrv}
\end{document}